\documentclass[a4paper,12pt]{amsart}
\usepackage{amsmath,latexsym,amscd}
\title[Bicanonical and adjoint linear systems]
{Bicanonical and adjoint linear systems on surfaces of general
type}
\author{Meng Chen and Eckart Viehweg}
\address{Institute of Mathematics, Fudan University,
Shanghai, 200433, PR China}
\email{mchen@fudan.edu.cn}

\address{Universit\"{a}t Essen, FB6 Mathematik, 45117 Essen, Germany}
\email{viehweg@uni-essen.de}

\thanks{This work has been supported by the ``DFG-Schwerpunktprogramm
Globale Methoden in der Komplexen Geometrie'' and by the DFG-NSFC
Chinese-German project ``Komplexe Geometrie''. The first named
author is supported by the National Natural Science Foundation of
China (Key Project No. 10131010) and by the Shanghai Scientific
$\&$ Technical Commission (Grant 01QA14042).}

\newcommand{\roundup}[1]{\ulcorner{#1}\urcorner}
\newcommand{\rounddown}[1]{\llcorner{#1}\lrcorner}

\newcommand{\sO}{{\mathcal O}}
\newcommand{\sE}{{\mathcal E}}

\newcommand{\Q}{{\mathbb Q}}
\newcommand{\C}{{\mathbb C}}
\newcommand{\bP}{{\mathbb P}}
\newtheorem{thm}{Theorem}[section]
\newtheorem{lem}[thm]{Lemma}

\newtheorem{prop}[thm]{Proposition}
\newtheorem{add}[thm]{Addendum}
\newtheorem{setup}[thm]{}
\theoremstyle{definition}

\newtheorem{question}[thm]{Question}
\newtheorem{exmp}[thm]{Example}

\newtheorem{rem}[thm]{Remark}
\newtheorem{claim}[thm]{Claim}
\numberwithin{equation}{thm}
\begin{document}
\begin{abstract}
This note contains a new proof of a theorem of Gang Xiao saying
that the bicanonical map of a surface $S$ of general type is
generically finite if and only if $p_2(S)>2$. Such properties are
also studied for adjoint linear systems $|K_S+L|$, where $L$ is
any divisor with $h^0(S,\sO_S(L))\geq 2$.
\end{abstract}
\maketitle
\section*{\bf Introduction}

Let $S$ be a complex minimal surface of general type. Since
$$K_S^2+1-q(S)+p_g(S) \geq 2$$
the Riemann-Roch Theorem implies that $p_2(S)\geq 2$. If
$p_2(S)=2$, the bicanonical map is composite with a pencil. It is
the aim of this note, to give an alternative proof of the Theorem
of G. Xiao, stating the converse.

\begin{thm}[Theorem 1 of \cite{X2}] \label{T:Xiao} Let $S$ be a
minimal projective surface of general type. Then the bicanonical
map of $S$ is generically finite if and only if $p_2(S)>2$.
\end{thm}

The proof of G. Xiao depends on his study of genus $2$ fibration
over curves and on Horikawa's classification of the possible
degenerations. We choose a different approach, and we will deduce
the Theorem from vanishing theorems for $\Q$-divisors, using in
addition just some well known and fundamental properties of
surfaces of general type.

We present such a new proof mainly as an interesting application
of the $\Q$-divisor method used for similar problems in higher dimensional
birational geometry (see for example \cite{MC}).
Using more involved results on surfaces, there are
other, slightly shorter proofs of Xiao's Theorem.

In the last section, we show that adjoint linear systems $|K_S+L|$
on surfaces of general type can only be composite with a pencil
of curves, if $L$ is a divisor with $h^0(S,\sO_S(L)) \leq 2$. We
discuss some examples, showing that this bound is sharp.
This result may be applied to study on 3-folds (see for example
\cite{MC}).\\

For a linear system $|L|$ on a surface $S$ the induced rational
map is denoted by $\varphi_L$. The linear system is composite
with a pencil of curves, if $\dim \varphi_L(S)=1$. The symbol
$\equiv$ stands for the numerical equivalence of divisors,
whereas $\sim$ denotes the linear equivalence. $K_S$ denotes the
canonical divisor, and if $f:S\to B$ is a surjective morphism,
$K_{S/B}=K_S - f^* K_B$. The base field is ${\mathbb C}$.

\section{\bf Proof of Theorem \ref{T:Xiao}}

Recall the Kawamata-Viehweg vanishing theorem (see \cite{Ka} or
\cite{V1}).

\begin{thm}[see \cite{E-V}, p. 49] \label{E-V}  Let $X$ be a smooth
projective variety and $L$ a divisor on $X$. Assume that $D$ is an
effective ${\mathbb Q}$-divisor with normal crossing supports
such that one of the following holds true:
\begin{enumerate}
\item[(i)] $L-D$ is nef and big.
\item[(ii)] $L-D$ is nef and $\kappa(L-\rounddown{D})=\dim X$.
\end{enumerate}
Then $H^i(X,\sO_S(K_X+L-\rounddown{D}))=0$ for all $i>0$.
\end{thm}

\begin{rem}
As well known, on surfaces, one may apply the vanishing theorems
without the assumption "normal crossings". In fact, if
$\tau:X'\to X$ is a blowing up, with $\tau^*D$ a normal crossing
divisor, then
$$R^i\tau_*\sO_{X'}(K_{X'}+\tau^*L-\rounddown{D'})=0,\mbox{ \ \ for }
i>0,$$ and for $i=0$ it coincides with
$\sO_X(K_X+L-\rounddown{D})$ in codimension one. If $X$ is a
surface, for $i>0$
\begin{gather*}
0=H^i(X',\sO_{X'}(K_{X'}+\tau^*L-\rounddown{D'}))\\
=H^i(X,\tau_*\sO_{X'}(K_{X'}+\tau^*L-\rounddown{D'}))
=H^i(X,\sO_X(K_X+L-\rounddown{D})).
\end{gather*}
\end{rem}

We will also use the following simple observation, due to Xiao
(see \cite{X2}, Lemme 8).

\begin{lem}\label{Xiao2}
Let $S$ be a minimal surface of general type with $q(S)=0$ and
$K_S^2\leq 2$. Let $\theta$ be a non-trivial invertible torsion
sheaf on $S$. Then $H^1(S,\theta)=0$.
\end{lem}
\begin{proof}
There exists an \'etale cover $\tau:T\to S$ with $\tau^*\theta =
\sO_T$, hence $\theta$ is a direct factor of $\tau_*\sO_T$. Since
$K_S^2 \leq 2 \leq 2 \chi(\sO_S)$ Corollary 5.8 of \cite{Bea}
implies that the fundamental group of $S$ is finite, hence the
one of $T$ as well. Then both $H^1(T,\sO_T)$ and $H^1(S,\theta)$
are zero.
\end{proof}

As a first step, let us reduce the proof of Theorem \ref{T:Xiao}
to the case $p_2(S)=3$.

\begin{prop}\label{P:2.1}
Let $S$ be a minimal smooth surface of general type.  Then
\begin{enumerate}
\item[(1)] the bicanonical map of $S$ is generically finite if
$p_2(S)\ge 4$;
\item[(2)] the linear system $|2K_S|$ is not
composite with an irrational pencil of curves for $p_2(S)=3$.
\end{enumerate}
\end{prop}
\begin{proof}
Suppose for some $S$ with $p_2(S) \ge 2$ the linear system
$|2K_S|$ is composite with a pencil, or for $p_2(S)=3$ with an
irrational pencil. Let $\pi: S'\longrightarrow S$ be any
birational modification such that $|2\pi^*(K_S)|$ defines a
morphism $\phi'_2$ and let $B_2'$ be its image. Consider the
Stein factorization
$$\phi'_2: S'\overset f\longrightarrow  B_2\longrightarrow  B_2'.$$
For some fibres $C_i$ of $f$ and for a general fibre $C$. We may
write
$$\pi^*(2K_S)\sim \sum_{i=1}^{a}C_i+Z_2 \equiv a\cdot C + Z_2,$$
where $Z_2$ is the fixed part. By assumption on the the smooth
curve $B_2$ the sheaf $f_*(\sO_{S'}(2K_{S'}))$ is invertible of
degree $a$ and the space of its global sections is of dimension
$\ge 4$, or of dimension $\ge 3$ if $B_2\neq \bP^1$. In both cases
one finds $a\ge 3$.

Set $G=\pi^*(K_S)-\frac{1}{a}Z_2.$ We have $K_{S'}+\roundup{G}\le
K_{S'}+\pi^*(K_S)$ and
$$G-C\equiv \frac{a-2}{a}\pi^*(K_S)$$
is nef and big. Thus \ref{E-V} implies that
$$|K_{S'}+\roundup{G}||_C=|K_C+D|,$$
for some divisor $D=\roundup{G}|_C$ of positive degree on the
curve $C$. The genus of $C$ can not be zero or one, hence $h^0(C,
K_C+D)\ge 2$. This implies that the morphism given by
$|K_{S'}+\pi^*(K_S)|$ can not factor through $f$, a contradiction.
\end{proof}

\begin{prop}\label{P:2.2} Let $S$ be a smooth minimal surface of general
type with $p_2(S)=3$. Assume that $|2K_S|$ is composite with a
pencil of curves. Then
\begin{enumerate}
\item[(i)] $K_S^2=2$ and $p_g(S)=q(S)\leq 1$.
\item[(ii)] $|2K_S|$ is composite with a rational pencil of curves of
genus $2$.
\item[(iii)] $|2K_S|$ defines a morphism on $S$, i.e. the movable part of
$|2K_S|$ is base point free.
\item[(iv)] Let $E$ be a component of the fixed part of $|2K_S|$. Then
$E\cdot K_S=0$ and $E$ is a $(-2)$ curve.
\end{enumerate}
\end{prop}
\begin{proof} Since $p_2(S)=3$ one has $p_g(S)\leq 2$.
The Riemann-Roch theorem and the positivity of the
Euler-Poincar\'e characteristic imply that
$$0 < K_S^2= 3 - 1 +q(S)-p_g(S) \leq 2.$$
By \cite{Bom}, Theorems 11 and 12, $q(S)=0$ if either $K_S^2=1$ or
if $K_S^2=p_g(S)=2$. Hence in order to prove (i), one just has to
exclude the case $K_S^2=1$, $p_g(S)=1$ and $q(S)=0$.

Since $p_2(S)=3$ Proposition \ref{P:2.1} implies that $|2K_S|$ is
composite with a rational pencil of curves. Let $\pi:S'\to S$ be
again a minimal birational modification such that $|2K_{S'}|$
defines a morphism $f:S'\to \bP^1$. The sheaf $f_*\sO_S(2K_S)$ is
invertible of degree two, hence we may write
$$2K_{S'}\sim 2C' + Z'_2$$
for a general fibre $C'$ of $f$. Set $C=\pi_*(C')$ and
$Z_2=\pi_*(Z'_2)$, then $2K_S\sim 2C+Z_2.$

If $K_S^2=1$ one has $C^2 \leq K_S\cdot C \leq 1$. Since the genus
of $C$ is at least two, $K_S\cdot C+C^2\ge 2$, which implies
$K_S\cdot C=C^2=1$ and $K_S^2\cdot C^2 =(K_S\cdot C)^2$. By the
Index Theorem $K_S\equiv C$. As shown in \cite{Bom} or
\cite{Catan} the condition $K_S^2=p_g(S)=1$ implies that on $S$
numerical equivalence coincides with linear equivalence. Hence
$K_S\sim C$, a contradiction since $p_g(S)\neq h^0(S,\sO_S(C))=2$.

Up to now, we obtained (i). For (iii) suppose that $\pi$ can not
be chosen to be an isomorphism, hence $C^2>0$. Then $2=K_S^2\ge
K_S\cdot C\ge C^2$. On the other hand, the index theorem gives
$$K_S^2\cdot C^2\le (K_S\cdot C)^2.$$
Since $K_S\cdot C + C^2$ is
even, one finds $K_S^2=K_S\cdot C=C^2=2$, hence $K_S\equiv C$, and
$Z_2=0$.

Assume $p_g(S)=1$. Let $D\in|K_S|$ be the unique effective
divisor. Then there are two fibers $C'_1$ and $C'_2$ of $f$ such
that, for $C_i=\pi(C'_i)$ one has $2D=C_1+C_2$.  If $C_1\ne C_2$,
then the $C_i$ are both 2-divisible for $i=1,\ 2$ and $D\equiv
2P$, where $P$ is a divisor. This implies $D^2\ge 4$, a
contradiction. If $C_1=C_2$, then $D=C_1$ and thus
$h^0(S,\sO_S(D))=2$, again a contradiction.

Assume $p_g(S)=0$, hence $q(S)=0$. Then the sheaf
$$\theta ={\mathcal O}_S(K_S-C)$$
is a non-trivial invertible torsion sheaf
on $S$. The Riemann-Roch Theorem implies $h^1(S,\theta )=1$,
contradicting Lemma \ref{Xiao2}.

So (iii) holds true and we may choose $S'=S$. Since for a general
fibre $C$ of $f$ one has $g(C) \ge 2$ and $K_S\cdot C \leq
K_S^2=2$, one finds $g(C)=2$, and $Z_2\cdot K_S = 0$.
\end{proof}

\begin{proof}[Proof of Theorem \ref{T:Xiao}.] By \ref{P:2.1} and \ref{P:2.2}
it remains to show, that there can not exist a surface with:
\begin{setup}\label{ass}
$S$ is a minimal surface of general type with $p_2(S)=3$, with
$K_S^2=2$ and with $p_g(S)=q(S)\leq 1$. The bicanonical map is a
genus two fibration $f:S \longrightarrow \bP^1$.
\end{setup}
Writing again $Z_2$ for the fixed part of $|2K_S|$ and $C$ for a
general fibre of $f$, one has $2K_S\sim 2C+Z_2$. Let $Z_v\leq
Z_2$ be the largest effective divisor contained in fibres of $f$,
and $Z_h=Z_2-Z_v$ the horizontal part of $Z_2$. In particular
$2C\cdot K_S= C\cdot Z_h = 4$. We will study step by step the
divisors $Z_v$ and $Z_h$.

\begin{claim}\label{C:1} The maximal multiplicity $a$ in $Z_2$ of an
irreducible component is two.
\end{claim}
\begin{proof} Suppose $a>2$, and denote by $\Gamma$ the total sum
of reduced components of multiplicity $a$ in $Z_2$. We may write
$$\Gamma=\Gamma_1+\cdots+\Gamma_s,$$
where the $\Gamma_i$ are connected pairwise disjoint.
\ref{P:2.2}, (iv), implies that each $\Gamma_i$ is a connected
tree of rational curves, thus 1-connected. We may replace $2C$ by
the sum of two different general fibres of $f$, say $C_1$ and
$C_2$. Then
$$K_S-\frac{1}{a}C_1-\frac{1}{a}C_2-\frac{1}{a}Z_2$$
is nef and big, and \ref{E-V} implies that
$$\textstyle H^1(2K_S-\Gamma_1-\cdots-\Gamma_s)=
H^1(2K_S+\roundup{-\frac{1}{a}C_1-\frac{1}{a}C_2-\frac{1}{a}Z_2})=0.$$
Thus we have a surjective map
$$H^0(S, 2K_S)\longrightarrow
H^0(\Gamma_1, {\mathcal O}_{\Gamma_1})\oplus\cdots\oplus
H^0(\Gamma_s, {\mathcal O}_{\Gamma_s})=\bigoplus^s\C,$$
contradicting $\Gamma\leq Z_2$.
\end{proof}

\begin{claim}\label{C:2} The horizontal part $Z_h$ of $Z_2$ is
either reduced, or $Z_h=2H$ for an irreducible (-2) curve $H$.
\end{claim}
\begin{proof}
If not, there is an irreducible curve $H_1$ with $Z_h-2H_1\neq
0$. By \ref{C:1} the multiplicities occurring in $Z_2$ are at
most $2$, and $Z_h\cdot C=4$ implies that either $Z_h-2H_1=2H_2$
for a reduced $(-2)$-curve $H_2$, or $Z_h-2H_1$ is reduced. Let
us write $H_2=0$ in the second case, such that in both cases
$$\textstyle \frac{1}{2}Z_h - \rounddown{\frac{1}{2}Z_h} +
H_2\neq 0. $$
Consider the effective $\Q$-divisor
$G=\frac{1}{2}(Z_2 - H_2)$. Obviously
$$K_S-G \equiv C+ \frac{1}{2}H_2$$
is nef. On the other hand,
$$\textstyle 2(K_S-\rounddown{G}) \ge 2C+Z_h - 2\rounddown{\frac{1}{2}Z_h} +
2H_2$$ is big. By the vanishing theorem \ref{E-V}, we have
$$H^1(S, 2K_S-\rounddown{G})=0.$$
The divisor $\rounddown{G}\geq H_1$ is again the sum over reduced
connected trees $\Gamma_i$ of $(-2)$-curves, say
$$\rounddown{G}= \Gamma_1+\cdots+\Gamma_s.$$
Thus we have a surjective map
$$H^0(S, 2K_S)\longrightarrow
H^0(\Gamma_1, {\mathcal O}_{\Gamma_1})\oplus\cdots\oplus
H^0(\Gamma_s, {\mathcal O}_{\Gamma_s})=\bigoplus^s\C,$$
contradicting $0 < 2\rounddown{G} \leq 2G \leq Z_2$.
\end{proof}

\begin{claim}\label{C:3}
$Z_h$ is either the sum of $4$ disjoint sections of $f$ or twice
an irreducible curve $H$. Moreover $Z_v=0$ in both cases.
\end{claim}
\begin{proof}
If $Z_h=2H$ for an irreducible curve $H$, one has $Z_h^2=-8$.
Otherwise \ref{C:2} only leaves the possibility $Z_h=H_1+ \cdots
+ H_t$, for $t\leq 4$. In this case, $Z_h^2 \geq -2t \geq -8$, and
$Z_h^2=-8$ if and only if $t=4$ and $H_i\cdot H_j=0$ for $i\neq
j$. The inequality
\begin{equation}\label{ineq1}
0=2K_S\cdot Z_h=8+Z_v\cdot Z_h+Z_h^2, \end{equation} implies
$Z_h^2\le -8$, and we obtain the first part of \ref{C:3}.

In both cases (\ref{ineq1}) is an equality, hence $Z_v\cdot
Z_h=0$. Finally the equality $$0=2K_S\cdot Z_v=2C\cdot
Z_v+Z_v^2+Z_v\cdot Z_h$$ implies $Z_v^2=0$ and by the Index
theorem $Z_v\equiv 0$. Since $Z_v\ge 0$ one finds $Z_v=0$.
\end{proof}
\begin{claim}\label{C:4}
In \ref{C:3} the case $Z_h=2H$ does not occur, and
$$Z_h=H_1+\cdots +H_4$$ implies $p_g(S)=q(S)=0$.
\end{claim}
\begin{proof}
Assume that $p_g(S) = 1$, and let $D$ denote the effective
canonical divisor. Then $2D=C_1+C_2+Z_h$ for fibres $C_i$ of $f$.
First of all this implies that the multiplicity of $Z_h$ is
divisible by $2$, hence $Z_h=2H$, and $C_1+C_2$ must be divisible
by $2$, as well. Since for any divisor $B$ the intersection number
$B^2+B\cdot K_S$ must be even, and since $C_i\cdot K_S = 2$ the
fibres $C_i$ can not be divisible by two. Hence $C_1=C_2$ and
$D=C_1+H$, a contradiction since $p_g(S) < h^0(S,\sO_S(D))=2$.

If $p_g(S)=0$, then by \ref{P:2.2} (i) $q(S)=0$. In case $Z_h=2H$
one finds $K_S\equiv C + H$ and $\theta = \sO_S(K_S-C-H)$ is a
$2$-torsion sheaf. The Riemann-Roch Theorem implies that
$h^1(S,\theta)=1$, contradicting \ref{Xiao2}.
\end{proof}
It remains to exclude the existence of a surface with:
\begin{setup}\label{ass2}
$S$ is a minimal surface of general type, $f:S \to \bP^1$ the
bicanonical map and for a fibre $C$ of $f$ and for pairwise
disjoint $(-2)$ curves $H_1,\ldots,H_4$
$$2K_{S/\bP^1}=6C+H_1+\cdots + H_4.$$
\end{setup}
Let us write $H=H_1+\cdots+H_4$. On some open dense subset
$U\subset \bP^1$ there is a natural involution $\iota$ on
$f^{-1}(U)$ with quotient $f^{-1}(U) \to \bP^1\times U$. Since
$S$ is minimal $\iota$ extends to an involution on $S$, denoted
again by $\iota$. The equality
$$0=2K_S\cdot \iota(H_i)=2C\cdot \iota(H_i)+(H_1+H_2+H_3+H_4)\cdot \iota(H_i)$$
implies that $\iota(H_i)\in \{H_1, H_2, H_3, H_4\},$ hence
$\iota(H)=H$. For $U$ small enough, each effective bicanonical
divisor of  $f^{-1}(U)$ is the pullback of a divisor on
$\bP^1\times U$, hence none of the $H_i$ can be fixed under
$\iota$. Renumbering we may assume that $\iota(H_1)=H_2$ and
$\iota(H_3)=H_4$.

Let $E$ be any $(-2)$-curve on $S$, not equal to one of the $H_i$.
The equality
$$0=2K_S\cdot E=2C\cdot E+(H_1+H_2+H_3+H_4)\cdot E$$
implies that $H_i\cdot E=0$ for all $i$. Hence $E$ is a component
of a fibre not meeting the $H_i$.

On the other hand let $E$ be any component of a fibre of $f$. If
$E$ does not meet $H$, then $E\cdot K_S=0$, hence $E$ is a
$(-2)$-curve.

The morphism $\delta : S \to S'$ to the relative minimal model
contracts exactly the $(-2)$ curves of the fibres. Hence all
fibres of ${f'}:S'\to \bP^1$ are reduced and all of their
components $E'$ meet $H'=\delta(H)$. Moreover the intersection
number $E\cdot K_S=E\cdot H$ on $S$ is even. So the reducible
fibres of ${f'}$ have at most two components $E'_1$ and $E'_2$, both
meeting $H'$ in two points. The components $E'_1$ and $E'_2$ need
not be Cartier divisors. However $E'_1+E'_2$ is Cartier, as well
as the images $H'_i$ of the $H_i$.

We write $\iota'$ for the automorphism of $S'$ induced by $\iota$.
Since $p_g(S)=q(S)=0$ the direct image $f_*\sO_S(K_{S/\bP^1})=
\sO_{\bP^1}(1)^{\oplus 2}$. Consider the restriction map
$$
\eta:{f'}_*\sO_{S'}(K_{S'/\bP^1}) = \sO_{\bP^1}(1)^{\oplus 2}
\longrightarrow \sO_{H'_1}(2)= \sO_{H_i}(K_{S'/\bP^1}\cdot H'_1).
$$
Since $\sO_C(K_C)$ is generated by global sections $\eta$ is
non-zero, hence its kernel is isomorphic to
$\sO_{\bP^1}(\epsilon)$, for $\epsilon = 0$ or $1$. Let $\sigma'$
be a general section of ${\rm Ker}(\eta)$, and let $\sigma$ be
the induced section of $\sO_{S'}(K_{S'/\bP^1})$. By construction
$H'_1$ lies in the zero-locus $B$ of $\sigma$. For some open
dense $U\subset \bP^1$ the divisor $B|_{{f'}^{-1}(U)}$ is invariant
under $\iota'$. Then the section $\sigma$ is zero on $H'_1+H'_2$.
Altogether we found an effective Cartier divisor $D'$ with
$$
\epsilon\cdot C + H'_1 + H'_2 + D' \sim K_{S'/\bP^1}.
$$
By construction $D'$ does not contain a whole fibre. So it is
concentrated in the reducible fibres of ${f'}$. Let
${f'}^{-1}(p)=E'_1 + E'_2$ be one of such fibres, and let $\alpha_1
\cdot E'_1 + \alpha_2 \cdot E'_2$ be the part of $D'$
concentrated in ${f'}^{-1}(p)$. Then one of the $\alpha_i$ must be
zero, say $\alpha_1$, hence $\alpha_2>0$.

The divisor ${\iota'}^*(\alpha_2 \cdot E'_2)$ is the part of
${\iota'}^*(D')$ lying in ${f'}^{-1}(p)$. If ${\iota'}^*(E'_2)=E'_1$
$$
\alpha_2 \cdot E'_2 -  {\iota'}^*(\alpha_2 \cdot E'_2)=
\alpha_2\cdot E'_2 - \alpha_2\cdot E'_1
$$
is the part concentrated in ${f'}^{-1}(p)$ of a divisor, linear
equivalent to zero. Then the same holds true for
$$
\alpha_2 \cdot \delta^*(E'_2) -  \alpha_2 \cdot \delta^* (E'_1).
$$
Obviously this is not possible, hence $E'_i$ is invariant under
$\iota'$.

We may assume that $E'_1\cap H'_1\neq \emptyset$. The component
$E'_1$ meets exactly one of the other $H'_i$, and being invariant
under $\iota'$, this can only be $H'_2$. Write $D=\delta^*(D')$
and $E_i$ for the proper transform of $E'_i$. If $D'$ contains
$E'_2$, it can not contain $E'_1$, hence $D$ does not contain
$E_1$. Since
$$
\epsilon\cdot C + H_1 + H_2 + D \sim K_{S/\bP^1}
$$
one finds $1=E_1\cdot K_{S/\bP^1}\geq E_1\cdot (H_1+H_2) = 2$,
obviously a contradiction. So $D'$ only contains components of
reducible fibres meeting $H'_1$ and $H'_2$ but neither $H'_3$ nor
$H'_4$. So $D\cdot H_3=0$ and
$$
H_3\cdot ( \epsilon\cdot C + H_1 + H_2 + D) = \epsilon < H_3\cdot
K_{S/\bP^1} = 2,
$$
a contradiction.
\end{proof}

\section{\bf Adjoint linear systems}

Let $S$ be a surface of general type, not necessary minimal, and
let $L$ be a divisor on $S$. There are few criteria known, which
imply that $\varphi_{K_S+L}$ is generically finite, though the
linear system $|K_S+L|$ quite well understood (see for instance
\cite{Rdr} and \cite{Catan2}).

By \cite{X3}, for a surface $S$ of general type with $q(S)\geq 3$
the map $\varphi_{K_S}$ is generically finite, hence the same
holds true for $\varphi_{K_S+L}$ whenever $L\geq 0$. We will
prove here

\begin{prop}\label{P:K+L} Let $S$ be a smooth projective surface of
general type and let $L$ be an effective divisor on $S$ with
$h^0(S,\sO_S( L))>2$. Then $\varphi_{K_S+L}$ is generically
finite.
\end{prop}

If $h^0(S,\sO_S( L)) = 2$ obviously $|L|$ is composite with a
pencil. The method used to prove \ref{P:K+L} will also
show:

\begin{add}\label{A:K+L} Assume in
\ref{P:K+L} that $h^0(S,\sO_S( L)) = 2$. Then $\varphi_{K_S+L}$
is generically finite, except possibly in one of the following
cases:
\begin{enumerate}
\item[(a)] $p_g(S)=0$ and $|L|$ is composite with a rational
pencil of hyperelliptic curves.
\item[(b)] $0<q(S)\leq 2$ and $|L|$ is composite with a rational pencil
of curves of genus $g = q(S)+1$.
\end{enumerate}
\end{add}

The next two examples shows that the exceptional cases
\ref{A:K+L}, (a) and (b), really occur.

\begin{exmp}\label{E:3.5}
In \cite{X2}, p. 46 - 49, one finds an example of a surface $S$ of
general type with $p_g(S)=q(S)=0$ and $K_S^2=2$, having a pencil
$f:S\to \bP^1$ of curves of genus $2$. If $C$ denotes a general
fibre, then
$$H^0(S,\sO_S(K_S+C))=H^0(C,\sO_C(K_C))=\C^{\oplus 2},$$
and $|K_S+C|$ is composite with a rational pencil of genus $2$
curves.
\end{exmp}

\begin{exmp}\label{E:3.6}
Let $C$ be a smooth curve of genus 2, and let $\theta$ be an
invertible $2$-torsion sheaf on $C$, with $\theta \neq \sO_C$.
For $T={\bP}^1\times C$ let $p_1:T\to {\bP}^1$ and $p_2:T\to C$
be the projections. For $a\geq 3$ consider
$$\delta=p_1^*(O(a))\otimes p_2^*(\theta).$$
Since $\delta^2\cong \sO_T(D)$ for a non-singular divisor $D$, one
obtains a smooth double cover $\pi:S\to T$ with
$$\pi_*\sO_S(K_S)=\sO_T(K_T)\oplus \sO_T(K_T)\otimes \delta.$$
It is easy to see that $S$ is a minimal surface of general type,
and that $|K_S|$ is composite with a pencil of curves of genus
$3$. In fact $\varphi_{K_S}$ coincides with $f=p_1\circ \pi$. For
a general fiber $C$ of $f$, choose $L=C$. Then $h^0(S,\sO_S(
L))=2$, but $|K_S+L|$ is composite with the same pencil as
$|K_S|$.

Note that $f$ is an isotrivial family of curves of genus $3$, that
$$f_*\sO_S(K_S)=\sO_{\bP^1}(a-2)\oplus \sO_{\bP^1}(-2)^{\oplus 2},$$
and that $q(S)=2$.
\end{exmp}

In Examples \ref{E:3.5} and \ref{E:3.6} the divisor $L$ is nef,
but not big.

\begin{question}\label{Q} Does there exists a minimal surfaces $S$
of general type and a nef and big divisor $L$ on $S$ with
$h^0(S,\sO_S( L))=2$, for which $|K_S+L|$ is composite with a
pencil of curves?
\end{question}

Such examples exist on surfaces $S$ of smaller Kodaira dimension,
or on surfaces $S$ of general type for $h^0(S,\sO_S( L))=1$:

\begin{exmp}\label{E:3.2}
Let $f:S\rightarrow \bP^1$ be a family of elliptic curves
admitting a section $G$, and with $S$ non-singular and projective.
For a general fibre $C$ of $f$ choose $L_m=mF+G$. Then $L_m$ is
nef and big, whenever $m>{\rm Max}\{0,-\frac{G^2}{2}\}$, and
$h^0(S,\sO_S(L_m))=m+1$. However $|K_S+L_m|$ is always composite
with a pencil.
\end{exmp}

\begin{exmp}\label{E:3.4}
Let $S$ be a minimal surface of general type with $K_S^2=1$ and
$p_g(S)=q(S)=0$. Denote by $L$ a divisor numerically equivalent
to $K_S$. Then $h^0(S,\sO_S( L)) \le 1$ and $h^0(S,
\sO_S(K_S+L))=2$. Thus $|K_S+L|$ is automatically composite with a
rational pencil of curves. One may refer to \cite{Rd} for a
classification of such pairs $(S, L)$.
\end{exmp}

\begin{proof}[Proof of \ref{P:K+L} (and of \ref{A:K+L})]
Replacing $S$ by a blowing up, we may assume that the moving part
of $L$ has no fixed points, hence that $\varphi_{L}$ is a
morphism.

Let us first consider the case that $|L|$ is composite with a
pencil of curves. Take the Stein factorization
\begin{equation}\label{stein}
g:S\overset{f}{\longrightarrow }B\overset{\rho}{\longrightarrow}
\bP(H^0(S,\sO_S(L))),
\end{equation}
so $f$ is a pencil of curves of genus
$g\geq 2$. As in the proof of \ref{P:2.1} one easily sees that
$h^0(S,\sO(L))>2$ implies that $L\geq C_1+C_2$ for two fibres $C_i$
of $f$. The same holds true for $h^0(S,\sO_S(L))=2$, if $\rho$
is not an isomorphism. In both cases we may as well assume that
$L=C_1+C_2$.

As explained in \cite{E-V}, 7.18, Koll\'ars vanishing theorem
implies that the locally free sheaf $f_{*}\sO_S(K_{S/B})$ is
numerically effective, and that $\sE=f_{*}\sO_S(K_{S}+C_1+C_2)$
is generated by global sections. Hence the tautological sheaf
$\sO_{\bP(\sE)}(1)$ on the projective bundle $\bP(\sE)$ is
globally generated.

If the genus $g(B)>0$, as a tensor product of a numerically
effective vector bundle with an invertible sheaf of positive
degree, $\sE$ is ample.

If $B\cong \bP^1$ the sheaf $\sE=f_{*}\sO_S(K_{S/B})$ is a direct
sum of line bundles of non-negative degree, say $\nu_1 \leq \nu_2
\leq \cdots \leq \nu_g$. If $q(S)=0$, by the Leray spectral
sequence $H^1(\bP^1,f_*\sO_S(K_S))=0$, hence $\nu_1>0$. If
$q(S)\neq 0$, one has $p_g(S)>0$, hence $\nu_g \geq 2$.

Altogether, in both cases the sheaf $\sO_{\bP(\sE)}(1)$ is
globally generated and big. $\varphi_{K_S+L}$ factors like
\begin{equation}\label{rel}
S \overset{\varphi} \longrightarrow \bP(\sE)
\overset{\varphi'}\longrightarrow \bP^M,
\end{equation}
where $\varphi$ is the relative canonical map and $\varphi'$ the
rational map induced by global sections of $\sO_{\bP(\sE)}(1)$.
Since the genus of the fibres of $f$ is at least two, $\varphi$
is generically finite. $\sO_{\bP(\sE)}(1)$, as well as its
restriction to the closure of the image of $\varphi$, are globally
generated and big, hence $\varphi_{K_S+L}$ is generically
finite.\\

Before finishing the proof of \ref{P:K+L} let us look to the case
$$h^0(S,\sO_S(L))=2, \mbox{ and } B\overset{\cong}{\longrightarrow} \bP^1$$
in (\ref{stein}). Here we may assume that $L=C$ for a general fibre of $f:S \to
\bP^1$. Write again $f_{*}\sO_S(K_{S/B})$ as a direct sum of line
bundles of non-negative degrees $\nu_1 \leq \nu_2 \leq \cdots
\leq \nu_g$. If $\varphi_{K_S+L}$ is composite with a pencil,
\cite{X3} implies that $q(S)<3$. Note that $\nu_i = 0$ for $i=1,
\ldots , q(S)$.

If $p_g(S)>0$, one also knows that $\nu_g\geq 2$. Hence if $g>
q(S)+1$, the sheaf $f_{*}\sO_S(K_{S}+C)$ contains a subbundle
$\sE$ of rank $\geq 2$ which is globally generated and non
trivial, i.e. not the direct sum of copies of $\sO_{\bP^1}$. For
this bundle consider again the maps (\ref{rel}). The first one,
$\varphi$, is fibrewise given by $\geq 2$ independent sections of
the canonical linear system, hence it is generically finite.
Since $\sO_{\bP(\sE)}(1)$, and its restriction to the image of
$\varphi$ are again generated by global sections and big,
$\varphi'\circ \varphi$ is generically finite and one obtains
\ref{A:K+L}, for $p_g(S)>0$.

If $p_g(S)=0$, hence $q(S)=0$, then $\nu_1 = \cdots = \nu_g =1$,
and $\sE=f_{*}\sO_S(K_{S}+C)$ is trivial. Then
$\bP(\sE)=\bP^1\times\bP^{g-1}$ and in (\ref{rel}) $\varphi$ is
generically finite, whereas $\varphi'$ is the projection to the
second factor. The restriction of $\varphi_{K_S+L}$ to a smooth
fibre $F$ coincides with $|K_F|$. So for $F$ non hyperelliptic,
the assumption that $|K_S+L|$ is composite with a pencil, implies
that all smooth fibres $F$ are isomorphic and that
$(\varphi_L,\varphi_{K_S+L})$ is a
birational map $S \to \bP^1\times F$, a contradiction.\\

To finish the proof of \ref{P:K+L} it remains to consider the
case that $\varphi_L$ is generically finite. If $p_g(S)>0$, the
linear system $|L|$ is a subsystem of $|K_S+L|$, hence the latter
can not be composite with a pencil of curves.

For $p_g(S)=q(S)=0$, blowing up $S$ if necessary, we assume that
both, $\varphi_{K_S+L}$ and $\varphi_{L}$ are morphism, hence
that the movable parts $M$ of $K_{S}+L$ and $L^0$ of $L$ have no
fixed points. Replacing $L$ by $L^0$ we may assume $L$ to be big
and globally generated.

Take the Stein factorization
$$\varphi_{K_S+L}: S\overset{h}{\longrightarrow }B\longrightarrow
{\bP}(H^0(S,\sO_S((K_S+L)-1))).$$ If $\varphi_{K_S+L}$ is not
generically finite, $h$ is a fibration onto a smooth curve $B$
with general fibre $C$. One may write $M\sim\sum_{i=1}^aC_i$ for
fibres $C_i$ of $h$ and for $a\ge h^0(S,\sO_S(K_S+L))-1$. Noting
that
$$h^0(S,\sO_S(K_{S}+L))=\frac{1}{2}L\cdot(K_{S}+L)+\chi({\mathcal O}_{S})
=\frac{1}{2}L\cdot(K_{S}+L)+1$$
one obtains the inequality
$$L\cdot(K_{S}+L)\ge L\cdot M \ge (\frac{1}{2}L\cdot(K_{S}+L))(L\cdot C), $$
hence $1\le L\cdot C\le 2$.

Consider next the natural map
$$H^0(S,\sO_S(L))\overset{\alpha}{\longrightarrow }W\subset
H^0(C,\sO_C(L|_C)),$$ with $W$ the image of $\alpha$. Because
$|L|$ is not composite with a pencil,
$$h^0(C,\sO_C(L|_C))\ge \dim_{\C}W\ge 2.$$
Noting that the genus $g(C)\ge 2$, one has
$h^0(C,\sO_C(\Gamma))\leq j$ whenever $\Gamma$ is a divisor with
$$1 \leq \deg(\Gamma)\leq j.$$ Hence
$$h^0(C,\sO_C(L|_C))=\dim_{\C}W=L\cdot C=2.$$
This implies that $h^0(S,\sO_S(L-C))\ge 1$ and $L-C\ge 0$. Since
$$|K_{S}+C||_C=|K_C|,$$ one finds $\dim\varphi_{K_{S}+L}(C)=1$,
contradicting the choice of $C$ as a fibre of $h$.
\end{proof}

\section*{\bf Acknowledgement} This note grew out of discussions
on the $\Q$-divisor method used in \cite{MC} for threefolds of
general type, during a visit of the first named author at the
University of Essen. He thanks the members of the Department of
Mathematics, in particular H\'el\`ene Esnault, for their
hospitality and encouragement. He also would like to thank
Fabrizio Catanese and Margarida Mendes Lopes for fruitful
discussions on adjoint linear systems on surfaces. After receiving
a first version of this note, Margarida Mendes Lopes gave another proof
of \ref{T:Xiao}, using some more advanced
techniques from surface theory.


\end{document}